\documentclass[11pt]{article}
\usepackage{amsfonts}
\usepackage{hyperref}
\usepackage{epsfig}
\usepackage{graphicx}

\newtheorem{theorem}{Theorem}[section]
\newtheorem{definition}[theorem]{Definition}
\newtheorem{proposition}[theorem]{Proposition}
\newtheorem{lemma}[theorem]{Lemma}
\newtheorem{claim}[theorem]{Claim}
\newtheorem{corollary}[theorem]{Corollary}

\newtheorem{question}[theorem]{Question}
\newtheorem{observation}[theorem]{Observation}

\newenvironment{proof}[1][Proof]{ \noindent \textbf{#1: }}{$\Box$
\bigskip}

\begin{document}

\title{On the Smallest Sets Blocking Simple Perfect Matchings
in a Convex Geometric Graph}

\author{
Chaya Keller and Micha A. Perles \\
Einstein Institute of Mathematics, Hebrew University\\
Jerusalem 91904, Israel\\
%{\tt nkeller@math.huji.ac.il}\\
}

\maketitle

\begin{abstract}

In this paper we present a complete characterization of the
smallest sets which block all the simple perfect matchings in a
complete convex geometric graph on $2m$ vertices. In particular,
we show that all these sets are caterpillar graphs with a special
structure, and that their total number is $m \cdot 2^{m-1}$.

\end{abstract}

\section{Introduction}

In this paper we consider geometric graphs (i.e., graphs whose
vertices are points in the plane, and whose edges are segments
connecting pairs of vertices), and in particular, convex geometric
graphs (i.e., geometric graphs whose vertices are in convex
position in the plane).
\begin{definition}
A simple perfect matching (SPM) in a geometric graph on $2m$
vertices is a set of $m$ pairwise disjoint edges (i.e., edges that
do not intersect, not even in an interior point).
\end{definition}
A natural Tur\'{a}n-type question (considered, e.g.,
in~\cite{Kupitz}) is: what is the maximal possible number of edges
in a geometric graph on $2m$ vertices with no simple perfect
matching?

\medskip

\noindent An equivalent way to state the question is to consider
sets which ``block'' all the SPMs:
\begin{definition}
A set of edges in a geometric graph $G$ is called a blocking set
if it intersects (i.e., contains an edge of) every SPM of the
graph.
\end{definition}
Using this formulation, the question above is equivalent to the
question:
\begin{question}
What is the minimal size (i.e., number of edges) of a blocking set
in a complete geometric graph on $2m$ vertices?
\end{question}
It appears that the answer depends on the position of the vertices
of the graph in the plane. It is easy to show that there always
exists a blocking set of size $2m-1$, and there exists a
configuration in the plane for which $2m-1$ is the minimal
possible size. On the other hand, in an unpublished
work~(\cite{Perles1}), Perles proved that for any placement of the
vertices, the size of a blocking set is at least $m$. This lower
bound is attained (among other cases) in the case of convex
geometric graphs (CGG.) Indeed, consider a complete convex
geometric graph on $2m$ vertices, denoted in the sequel $CK(2m)$.
The vertices of the graph form a convex polygon. It is easy to see
that any set of $m$ consecutive edges on the boundary of the
polygon is a blocking set. The set of all edges of odd order
emanating from a single vertex is also clearly a blocking set of
size $m$.\footnote{For the sake of clarity, we present proofs of
these two straightforward claims in
Section~\ref{sec:sub:observations}.}

\bigskip

\begin{figure}[tb]
\begin{center}
\scalebox{0.6}{
\includegraphics{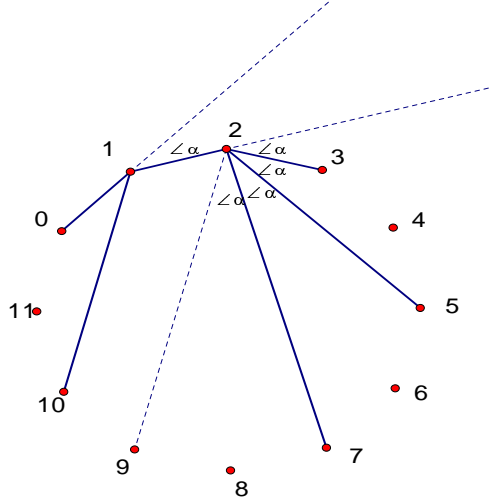}
} \caption{A blocker in $CK(12)$ with spine of length $t=3$. The
edges of the blocker are depicted by full (not dotted) lines. In
the notation of Theorem~\ref{Thm:Main}, $\epsilon_4=1$,
$\epsilon_5=2$, and $\epsilon_6=4$. The angle $\alpha$ is $\pi/6$
radians. The diagonal $[2,9]$ is parallel to the diagonal
$[1,10]$, and helps to depict the angle between the diagonals
$[2,7]$ and $[1,10]$.} \label{fig:T1}
\end{center}
\end{figure}

\noindent In this paper we present a complete characterization of
the blocking sets of size $m$ in $CK(2m)$, called in the sequel
{\it blockers}. It turns out that all these blockers are simple
subtrees of a special structure, called {\it caterpillars} (see,
e.g.,~\cite{Caterpillar1}).
\begin{definition}
A tree $T$ is a caterpillar (or a fishbone) if the derived graph
$T'$ (i.e., the graph obtained from $T$ by removing all leaves and
their incident edges) is a path (or is empty). A geometric
caterpillar is simple if it does not contain a pair of crossing
edges. A longest (simple) path in a caterpillar $T$ is called a
spine of $T$.
\end{definition}
Our main result is the following:
\begin{theorem}
Let $V$ be the set of vertices of a convex $2m$-gon $P$, labelled
cyclically from $0$ to $2m-1$, and let $G$ be the complete convex
geometric graph on $V$. Any blocker of $G$ is a simple caterpillar
graph whose spine lies on the boundary of the polygon and is of
length $t \geq 2$. If the spine ``starts'' with the vertex $0$ and
the edge $[0,1]$, then the edges of the blocker are:
\begin{equation}
\{[i-1,i]:1 \leq i \leq t\} \cup
\{[t+j-1-\epsilon_{t+j},t+j+\epsilon_{t+j}]:1 \leq j \leq m-t\},
\end{equation}
where the $\epsilon_i$ are natural numbers satisfying $1 \leq
\epsilon_{t+1} < \epsilon_{t+2} < \ldots < \epsilon_m \leq m-2$.

Conversely, any set of $m$ edges of the described form is a
blocker in $G$. \label{Thm:Main}
\end{theorem}
If the polygon is regular, then the direction of each consecutive
edge of the blocker, as listed above, is obtained from the
direction of the preceding edge by rotation by $\pi/m$ radians. In
the first $t$ edges, the ``back'' endpoint of each edge is the
``front'' endpoint of the previous edge. Starting with the
$t+1$-st edge, the ``back'' endpoint goes ``back'' (as reflected
by subtraction of the corresponding $\epsilon_i$), and the length
of the edge changes accordingly. An example of a blocker in
$CK(12)$ is presented in Figure~\ref{fig:T1}.

\bigskip

\noindent The proof of the theorem involves various techniques,
including examination of several specific classes of SPMs, as well
as inductive arguments.

\medskip

\noindent As an easy corollary of the structure theorem, we
enumerate the blockers in $CK(2m)$:
\begin{proposition}
Let $G=CK(2m)$ be a complete convex geometric graph on $2m$
vertices. The number of blocking sets of size $m$ in $G$ is $m
\cdot 2^{m-1}$.
\end{proposition}

%\noindent Finally, we consider blockers with more than $m$ edges,
%which are minimal with respect to inclusion. We prove that for any
%convex configuration of the vertices, there exist such blockers of
%size $m^2/4$. However, we do not know whether this size is the
%maximal possible.

\section{Preliminaries}
\label{sec:Preliminaries}

In this section we introduce several basic definitions and
observations, and consider two specific classes of SPMs that will
be used in the proof of Theorem~\ref{Thm:Main}.

\subsection{Definitions and Observations}
\label{sec:sub:observations}

\begin{definition}
Let $V$ be the set of vertices of a convex $2m$-gon $P$, $m \geq
2$, labelled cyclically from $0$ to $2m-1$.
\begin{itemize}
\item A \textbf{half-boundary} of $P$ is a set of $m$ consecutive
boundary edges.

\item The \textbf{order} of an edge $[i,i+k]$ (where the addition
is modulo $2m$) is $\min(k,2m-k)$. The boundary edges of $P$ are,
of course, of order $1$. We call the non-boundary edges, i.e., the
edges that are diagonals of $P$, \textbf{interior edges}.

\item Let $[i-1,i]$ and $[i+k-1,i+k]$ be two boundary edges of $P$
(where $0<k<2m$, and the addition within the edges is modulo
$2m$). The \textbf{(directed) distance} from $[i-1,i]$ to
$[i+k-1,i+k]$ is $k$. In particular, the distance from a boundary
edge to its immediate successor is 1. (Note that the distance from
$[i+k-1,i+k]$ to $[i-1,i]$ is $2m-k$.)
\end{itemize}
\end{definition}

\begin{observation} \label{obs:basic}
Let $V$ be the set of vertices of a convex $2m$-gon $P$, labelled
cyclically from $0$ to $2m-1$, and let $G$ be the complete convex
geometric graph on $V$. Then:
\begin{enumerate}
\item Any blocker in $G$ contains at least two boundary edges.

\item The set of all edges of odd order emanating from a single
vertex is a blocker in $G$.

\item Any set of $m$ consecutive boundary edges of $P$ is a
blocker in $G$.

\end{enumerate}
\end{observation}

\begin{proof}
\begin{enumerate}
\item The boundary of $P$ is the disjoint union of two SPMs:
$\{[2i,2i+1]:i=0,1,\ldots,m-1\}$, and
$\{[2i+1,2i+2]:i=0,1,\ldots,m-1\}$. In order to intersect these
two SPMs, any blocker has to contain at least two boundary edges.

\item Note that all the edges in any SPM are of odd order. Indeed,
an edge $[i,j]$ of an SPM $M$ divides the remaining vertices of
$P$ into sets $V_1,V_2$ of sizes $j-i-1$ and $2m-2-(j-i-1)$. Since
the edges of $M$ do not intersect, the two vertices of any other
edge are in the same set (either both in $V_1$ or both in $V_2$).
As $M$ ``covers'' all the vertices of $P$, it follows that each of
the sets $V_1,V_2$ contains an even number of vertices. Hence,
$j-i-1$ is even, and thus the order of the edge $[i,j]$ is odd.

Let $B$ be the set of odd-order edges emanating from the vertex
$v$, and let $M$ be an SPM. Since $M$ is a prefect matching, it
contains an edge emanating from $v$. By the explanation above,
this edge is of odd order, so it is included in $B$. Thus, $B$
intersects $M$, as asserted.

\item Assume w.l.o.g. that the set is
$B=\{[0,1],[1,2],\ldots,[m-1,m]\}$, and let $M$ be an SPM. By the
pigeonhole principle, $M$ contains an edge with both vertices in
$\{0,1,2,\ldots,m\}$. Let $[i_0,j_0]$ (for $i_0<j_0$) be a
``shortest'' (i.e., having the smallest order) edge with this
property. Since the edges of $M$ cover all the vertices and do not
intersect, each of the vertices in the set
$\{i_0+1,i_0+2,\ldots,j_0-1\}$ is ``connected'' by $M$ to another
vertex in this set. However, an edge that connects two such
vertices is shorter than $[i_0,j_0]$, contradicting the assumption
above. Thus, the set $\{i_0+1,i_0+2,\ldots,j_0-1\}$ is empty, so
$[i_0,j_0]$ is a boundary edge, which is contained in $B$.
Therefore $B$ intersects any SPM, as asserted.
\end{enumerate}
\end{proof}

\subsection{Parallel SPMs}
\label{sec:sub:parallel}

We start with a combinatorial generalization of the notion of
parallel edges.

\medskip

\noindent If the polygon $P$ (that consists of the vertices and
boundary edges of $CK(2m)$) is regular, then its edges and
diagonals have $2m$ directions: $m$ directions of the boundary
edges and the diagonals of odd order, and $m$ directions of the
diagonals of even order. The directions define an equivalence
relation, whose equivalence classes consist of all the boundary
edges and diagonals of the same direction. The equivalence classes
of the first type (odd order) contain two boundary edges and $m-2$
diagonals, and the equivalence classes of the second type (even
order) contain $m-1$ diagonals. This equivalence relation can be
defined in a combinatorial way, that extends naturally to the
edges and diagonals of any convex polygon of even order.
\begin{definition}
Let $[p,q]$ and $[p',q']$ be disjoint segments connecting four
different vertices of a convex polygon $P$ on $2m$ vertices, such
that the order of the vertices on the boundary of the polygon is
$p,q,p',q'$. The segments are called ``parallel'' if the number of
boundary edges in the arc $\langle q,p' \rangle$ is equal to the
number of boundary edges in the arc $\langle q',p \rangle$.
\end{definition}
A special class of SPMs we consider consists of full equivalence
classes of the relation defined above.
\begin{definition}
The set of all edges which are parallel to a given boundary edge
is called a ``parallel SPM''. The parallel SPMs are of the form
$M_l=\{[i,j]:i+j \equiv 2l-1 (\bmod 2m)\}$, for all $1 \leq l \leq
m$.
\end{definition}
The sets $\{M_l\}_{l=1}^m$ are pairwise disjoint. Since a blocker
has only $m$ edges and intersects each of the parallel SPMs (i.e.,
each of the sets $M_l$), it must intersect each of the $M_l$-s in
exactly one edge. We thus get the following:
\begin{observation}
\label{obs:parallel}
Any blocker contains exactly one edge of each
of the equivalence classes of odd order.
\end{observation}

\subsection{Triangular SPMs}
\label{sec:sub:Triangular-SPM}

\begin{figure}[tb]
\begin{center}
\scalebox{0.5}{
\includegraphics{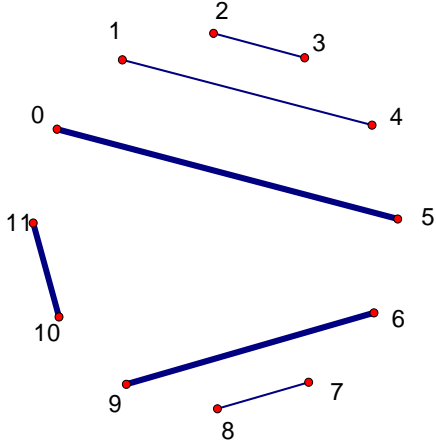}
} \caption{A triangular SPM in $CK(12)$, corresponding to the case
$i_0=0$, $a=3$, $b=2$, $c=1$ (in the above notations). The initial
(innermost) diagonals are drawn thick. The boundary edges of the
SPM are $[2,3],[7,8],$ and $[10,11]$, and the distances between
them are $5,3,$ and $4$, respectively.} \label{fig:T2}
\end{center}
\end{figure}

For any triple of positive integers $(a,b,c)$ with $a+b+c=m$ and a
``starting point'' $i_0$, $0 \leq i_0 \leq 2m-1$, consider the
triple of segments
\[
\Big([i_0,i_0+2a-1],[i_0+2a,i_0+2a+2b-1],[i_0+2a+2b,i_0-1]\Big),
\]
where the additions are taken modulo $2m$. Note that the segments
are pairwise disjoint diagonals (or edges) of the polygon $P$.
This triple of segments can be extended to an SPM by adding the
following segments:
\[
[i_0+\epsilon,i_0+2a-1-\epsilon], \qquad \epsilon=1,2,\ldots,a-1,
\]
\[
[i_0+2a+\epsilon,i_0+2a+2b-1-\epsilon], \qquad
\epsilon=1,2,\ldots,b-1,
\]
\[
[i_0+2a+2b+\epsilon,i_0-1-\epsilon], \qquad
\epsilon=1,2,\ldots,c-1.
\]
An SPM of this form is called a {\it triangular SPM} (see Figure~\ref{fig:T2}).

The boundary
edges of this triangular SPM are
\[
[i_0+a-1,i_0+a],[i_0+2a+b-1,i_0+2a+b],[i_0+2a+2b+c-1,i_0+2a+2b+c].
\]
The distances between these boundary edges, in cyclical order, are
$a+b,b+c,$ and $c+a$, and by assumption, all of them are less than
$m$. In the following proposition, that will be used in the proof
of our main theorem, we claim that the converse holds as well:
\begin{proposition}
For any triple of boundary edges
$\Big([i_1-1,i_1],[i_2-1,i_2],[i_3-1,i_3]\Big)$, $1 \leq
i_1<i_2<i_3<2m$, such that the distance from each one to the next
(in cyclical order) is less than $m$, there exists a triangular
SPM whose boundary edges are
$[i_1-1,i_1],[i_2-1,i_2],[i_3-1,i_3]$. \label{Prop:Triangular}
\end{proposition}
\begin{proof}
Denote the distances from each edge to the next, in cyclical
order, by $p,q,r$. That is,
\[
p=i_2-i_1, \qquad q=i_3-i_2, \qquad r=i_1+2m-i_3.
\]
By assumption, $0 < p,q,r < m$. Consider a set of edges that
consists of $a$ consecutive edges parallel to $[i_1-1,i_1]$ (i.e.,
$\{[i_1-1-\epsilon,i_1+\epsilon]\}$, $\epsilon =0,1,\ldots,a-1$),
$b$ consecutive edges parallel to $[i_2-1,i_2]$, and $c$
consecutive edges parallel to $[i_3-1,i_3]$. It is easy to see
that this set is an SPM if the following three equalities hold:
\[
(1) \mbox{  } a+b=p, \qquad (2) \mbox{  } b+c=q, \qquad (3) \mbox{
} c+a=r.
\]
Summing the equations we get $2(a+b+c)=p+q+r=2m$. Subtracting
equations (1),(2),(3) from the equation $a+b+c=m$, we get the
solutions $(a=m-q,b=m-r,c=m-p)$, and these are indeed positive
integers. Thus, the set of edges defined above with
$a=m-q,b=m-r,c=m-p$ is an SPM whose boundary edges are
$[i_1-1,i_1],[i_2-1,i_2],[i_3-1,i_3]$, as claimed.
\end{proof}

\begin{figure}[tb]
\begin{center}
\scalebox{0.6}{
\includegraphics{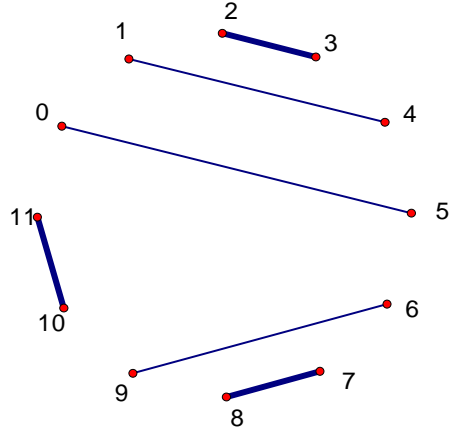}
} \caption{Construction of an SPM in $CK(12)$ given three boundary
edges. In this figure, $i_1=3$, $i_2=8$, and $i_3=11$. The given
boundary edges are drawn thick. The obtained values are $a=3,b=2,$
and $c=1$, and the obtained SPM is the same as in
Figure~\ref{fig:T2}.} \label{fig:T3}
\end{center}
\end{figure}

The construction of a triangular SPM from three given boundary
edges is exemplified in Figure~\ref{fig:T3}.

\section{Proof of Theorem~\ref{Thm:Main}}
\label{sec:proof}

In this section we present the proof of our main theorem. We start
with an outline of the proof.

\subsection{Proof Outline}
\label{sec:sub:outline}

The key observation is that a characterization of the possible
{\it boundary edges} in a blocker leads to a full characterization
of the blockers. The main step in the proof is the following
lemma, characterizing the boundary edges of a blocker:
\begin{lemma}
The boundary edges of a blocker form a path of length $t$ on the
boundary of the polygon, $2 \leq t \leq m$. \label{Lemma:Main}
\end{lemma}
Lemma~\ref{Lemma:Main}, in turn, is proved in two steps. First we
prove:
\begin{lemma}
The boundary edges of a blocker are included in a half-boundary.
\label{Lemma:Half-Boundary}
\end{lemma}
We prove Lemma~\ref{Lemma:Half-Boundary} by showing that if the
boundary edges are not included in a half-boundary then there
exists an SPM of one of the two special kinds mentioned above
(``parallel'' and ``triangular'') that misses the blocker. Then we
deduce Lemma~\ref{Lemma:Main} from Lemma~\ref{Lemma:Half-Boundary}
by an inductive argument. Using Lemma~\ref{Lemma:Main}, we show
that if a set of $m$ edges is not a caterpillar with the specified
properties, then there exists an SPM that misses it. This proves
one direction of Theorem~\ref{Thm:Main}.

\medskip

\noindent The other direction of the theorem (asserting that any
caterpillar with the specified properties is a blocker) is proved
by double induction: A primary induction on $m$, and a secondary
(backward) induction on the number of boundary edges in the
caterpillar.

\subsection{Proof of Lemma~\ref{Lemma:Half-Boundary}}
\label{sec:sub:half-boundary}

We use the following technical lemma:
\begin{proposition}
Let S=$\{[i_1,i_1+1],[i_2,i_2+1],\ldots,[i_k,i_{k}+1]\}$ be a set
of $k$ boundary edges of $CK(2m)$, where $0 \leq
i_1<i_2<\ldots<i_k \leq 2m-1$. Then at least one of the following
holds:
\begin{enumerate}
\item $S$ contains two opposite edges (i.e., $i_{\nu}=i_{\mu}+m$
for some $\mu,\nu$, $1 \leq \mu < \nu \leq k$).

\item $S$ is included in a half-boundary (i.e., there exists a
$\mu$, $1 \leq \mu < k$, such that $i_{\mu}+m < i_{\mu+1}$, or
$i_k<i_1+m$).

\item $S$ contains three edges such that the distance from each
one to the next (in cyclical order) is less than $m$ (i.e., there
exist $1 \leq \mu < \nu < \tau \leq k$ such that $i_{\nu} <
i_{\mu} + m$ , $i_{\tau} < i_{\nu}+m$, and $i_{\mu}+2m <
i_{\tau}+m$).
\end{enumerate}
\label{Prop:Three-Options}
\end{proposition}
\begin{proof}
For $k=1,2$, it is easy to see that either (1) or (2) holds. Let
$k>2$, and assume that both (1) and (2) fail. Define $\mu_0=\max
\{\mu: 2 \leq \mu \leq k, i_{\mu}<i_1+m \}$. Since (2) fails, we
have $2 \leq \mu_0 <k$. Consider the edges
$[i_1,i_1+1],[i_{\mu_0},i_{\mu_0}+1],[i_{\mu_0+1},i_{\mu_0+1}+1]$.
These edges satisfy the requirements of (3). Indeed, by the
definition of $\mu_0$, $i_{\mu_0}<i_1+m$. Furthermore, we have
$i_{\mu_0+1} < i_{\mu_0} +m$, since otherwise either (1) or (2)
are satisfied. Finally, by the definition of $\mu_0$, $i_{\mu_0+1}
\geq i_1 +m$, and equality cannot hold since, by assumption, (1)
fails. Hence, $i_{\mu_0+1}>i_1+m$, and thus (3) holds. This
completes the proof.
\end{proof}

\medskip

\noindent Now we are ready to prove
Lemma~\ref{Lemma:Half-Boundary}. The formal statement of the lemma
is the following:
\begin{lemma}
Let
$B=\{[i_1,i_1+1],\ldots,[i_k,i_k+1],[i_{k+1},j_{k+1}],\ldots,[i_m,j_m]\}$
be a blocker in $CK(2m)$, where for all $1 \leq \mu \leq k$,
$[i_{\mu},i_{\mu}+1]$ is a boundary edge, and for all $k< \mu \leq
m$, $[i_{\mu},j_{\mu}]$ is not a boundary edge. Then the edges
$[i_1,i_1+1],\ldots,[i_k,i_k+1]$ are included in a half-boundary.
\label{Lemma':Half-Boundary}
\end{lemma}
\begin{proof}
Consider the set $E=\{[i_1,i_1+1],\ldots,[i_k,i_k+1]\}$ of
boundary edges of the blocker. Assume, w.l.o.g., that $0 \leq
i_1<i_2<\ldots<i_k \leq 2m-1$. By
Proposition~\ref{Prop:Three-Options}, at least one of the
following holds:
\begin{enumerate}
\item $E$ contains two opposite edges (i.e., there exist
$\mu,\nu$, $1 \leq \mu< \nu \leq k$, such that
$i_{\nu}=i_{\mu}+m$).

\item $E$ contains three edges such that the distance from each to
the next (in cyclical order) is less than $m$ (i.e., there exist
$1 \leq \mu < \nu < \tau \leq k$ such that $i_{\nu} < i_{\mu} + m$
, $i_{\tau} < i_{\nu}+m$, and $i_{\mu}+m < i_{\tau}$).

\item $E$ is included in a half-boundary (i.e., there exists a
$\mu$, $1 \leq \mu < k$, such that $i_{\mu}+m < i_{\mu+1}$, or
$i_k<i_1+m$).
\end{enumerate}
(1) is impossible, since by Observation~\ref{obs:parallel}, $B$
does not contain two parallel edges.

\medskip

\noindent If $E$ contains three edges $[i_{\mu},i_{\mu}+1],
[i_{\nu},i_{\nu}+1], [i_{\tau},i_{\tau}+1]$ such that the distance
from each one to the next (in cyclical order) is less than $m$,
then the triple of opposite edges,
\[
T= \Big( [i_{\mu}+m,i_{\mu}+1+m], [i_{\nu}+m,i_{\nu}+1+m],
[i_{\tau}+m,i_{\tau}+1+m] \Big),
\]
also has this property. By Proposition~\ref{Prop:Triangular}, the
triple $T$ can be extended to a triangular SPM $\tilde{T}$. Each
edge of $\tilde{T}$ is parallel to an edge in $T$, and the only
boundary edges of $\tilde{T}$ are the three edges of $T$. It
follows that our blocker $B$ misses the SPM $\tilde{T}$ entirely:
the only edges of $B$ that are parallel to edges in $T$ are the
boundary edges $[i_{\mu},i_{\mu}+1], [i_{\nu},i_{\nu}+1],$ and $
[i_{\tau},i_{\tau}+1]$, which are not in $\tilde{T}$. Hence, (2)
is also impossible. Therefore, we are left with (3), i.e., $E$ is
included in a half-boundary, as claimed.
\end{proof}

\subsection{Proof of Lemma~\ref{Lemma:Main}}

In the proof of the lemma we use the following inductive
technique.

Let $B$ be a blocker in $CK(2m)$. Consider a pair of consecutive
boundary edges $e,f$, such that $e \in B$ and $f \not \in B$. Such
a choice is possible, since by Observation~\ref{obs:basic}, the
number of boundary edges of $B$ is between 2 and $m$. Assume,
w.l.o.g., that $e=[2m-3,2m-2]$, and $f=[2m-2,2m-1]$.

\begin{figure}[tb]
\begin{center}
\scalebox{0.6}{
\includegraphics{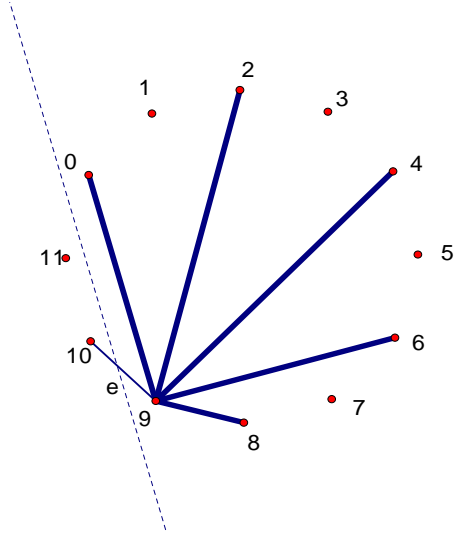}
} \caption{An illustration of the proof of
Claim~\ref{Claim:Inductive}. In this figure, $m=6$, $e=[9,10]$,
and $f=[10,11]$. The vertices of the induced subgraph $CK(10)$ are
$0,1,\ldots,9$, the edges of $B'$ are drawn thick, and $B=B' \cup
\{e\}$.} \label{fig:T4}
\end{center}
\end{figure}

Denote by $CK(2m-2)$ the geometric subgraph of $CK(2m)$ obtained
by omitting the endpoints of $f$. The boundary of $CK(2m-2)$ is
thus $\langle 0,1,2,\ldots,2m-3,0 \rangle$.
\begin{claim}\label{Claim:Inductive}
The set $B \setminus \{e\}$ is a blocking set of size $m-1$ (i.e.,
a blocker) in $CK(2m-2)$.
\end{claim}
\begin{proof}
Let $B'=B \cap E(CK(2m-2))$. $B'$ is obtained from $B$ by omitting
$e$ and any other edge that uses one of the vertices $2m-2,2m-1$
(see Figure~\ref{fig:T4}). If $B'$ is not a blocking set in
$CK(2m-2)$ then there exists an SPM in $CK(2m-2)$ that misses
$B'$, and thus misses $B$. Adding the edge $f$ to that SPM yields
an SPM in $CK(2m)$ that misses $B$, contradicting the assumption
that $B$ is a blocker. Hence, $B'$ is a blocking set in
$CK(2m-2)$. Clearly, $B' \subseteq B \setminus \{e\}$. This
inclusion must be an equality, i.e., $B'=B \setminus \{e\}$, since
a blocking set in $CK(2m-2)$ cannot have fewer than $m-1$ edges.
\end{proof}

The same argument yields immediately the following corollary.
\begin{corollary}
In the notations above, $B$ does not contain any edge that uses
one of the vertices of $f$ (i.e., $2m-2$ and $2m-1$), except $e$.
\end{corollary}

\begin{figure}[tb]
\begin{center}
\scalebox{0.6}{
\includegraphics{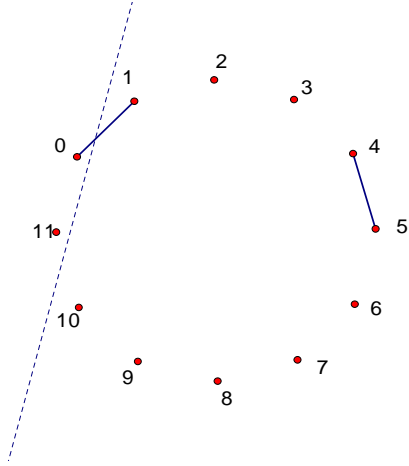}
} \caption{An illustration of Case~4 of the proof of
Lemma~\ref{Lemma':Main}. In this figure, $m=5$, and $i_t=4$. In
the graph $CK(12)=CK(2m+2)$, the boundary edges of the blocker are
$[0,1]$ and $[4,5]$. In the induced subgraph $CK(10)$ (to the
right of the dotted line), the boundary edges of the induced
blocker are $[4,5]$ and possibly also $[1,10]$, and both options
lead to a contradiction.} \label{fig:T7}
\end{center}
\end{figure}

Now we are ready to prove Lemma~\ref{Lemma:Main}. The formal
statement is the following:
\begin{lemma} \label{Lemma':Main}
Let $B$ be a blocker in $G=CK(2m)$. The boundary edges of $B$ are
consecutive (i.e., if $B$ contains $t$ boundary edges
$[i_1,i_1+1],\ldots,[i_t,i_{t}+1]$, then these edges can be
arranged in such a way that $i_{\mu+1}=i_{\mu}+1 (\bmod 2m)$, for
all $1 \leq \mu < t$).
\end{lemma}
\begin{proof}
The proof is by induction on $m$.

\textbf{The case $m=2$:} Any blocker in $CK(4)$ contains exactly
one edge of each $M_l=\{[i,j]:i+j=2l-1 (\bmod 4)\}$ ($l=1,2$). Thus, any
such blocker consists of two consecutive boundary edges.

\textbf{The case $m=3$:} If a blocker in $CK(6)$ contains three
boundary edges, then, by Lemma~\ref{Lemma:Half-Boundary}, these
edges are contained in a half-boundary, and thus are consecutive.
Since by Observation~\ref{obs:basic}, any blocker contains at
least two boundary edges, we are left with the case where the
blocker contains exactly two boundary edges. Assume on the
contrary that these edges are not consecutive. Since these edges
are not parallel either, we may assume, w.l.o.g., that they are
$[0,1]$ and $[4,5]$ (see Figure~\ref{fig:T6}). However, in this
case the blocker misses the SPM $([1,2],[3,4],[5,0])$.

\begin{figure}[tb]
\begin{center}
\scalebox{0.6}{
\includegraphics{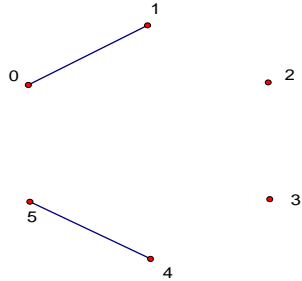}
}
\caption{An attempt to construct the boundary edges of a blocker in CK(6).}
\label{fig:T6}
\end{center}
\end{figure}

\textbf{The inductive step:} We assume that the claim holds for
any blocker in $CK(2m)$, and prove it for a blocker $B$ in
$CK(2m+2)$. Denote the set of boundary edges of $B$ by $E$,
$E=\{[i_1,i_1+1],\ldots,[i_t,i_t+1]\}$, and assume w.l.o.g. that
$[i_1,i_1+1]=[0,1]$, that the sequence $\{i_1,i_2,\ldots,i_t\}$ is
monotone increasing, and that $E \subset
\{[0,1],[1,2],\ldots,[m,m+1]\}$. (The last assumption is valid
since, by Lemma~\ref{Lemma:Half-Boundary}, $E$ is included in a
half-boundary). We perform a case-by-case analysis.

\textbf{Case 1: $[1,2] \in E$.} Since $[2m+1,0] \not \in E$ and
$[0,1] \in E$, by the inductive technique presented above
(Claim~\ref{Claim:Inductive}), the set $B \setminus \{[0,1]\}$ is
a blocker in the subgraph spanned by $\{1,2,\ldots,2m\}$. Hence,
by the inductive assumption, the boundary edges of $B \setminus
\{[0,1]\}$ form a connected set that contains the edge $[1,2]$,
and therefore $E$, the set of boundary edges of $B$, is also
connected. (Note that since, by assumption, $[2m-1,2m] \not \in
B$, the boundary edges of $B \setminus \{[0,1]\}$ in the induced
subgraph are $[1,2],[2,3],\ldots,[t-1,t]$ for some $t \leq m+1$,
and possibly also $[2m,1]$, but not $[2m-1,2m]$. Since $[2m,1]$ is
not a boundary edge in the original graph, it follows that
$E=\{[0,1],[1,2],[2,3],\ldots,[t-1,t]\}$, which is indeed a set of
consecutive edges.)

\textbf{Case 2: $[i_t-1,i_t] \in E$.} In this case we can use an
argument symmetric to that used in Case~1.

\textbf{Case 3: $[1,2],[i_t-1,i_t] \not \in E$, and $|E| \geq 3$.}
Let $[j,j+1] \in E \cap \{[2,3],[3,4],\ldots,[i_t-2,i_t-1]\}$. As
in Case~1, by the inductive technique presented above, $B
\setminus \{[0,1]\}$ is a blocker in the subgraph spanned by
$\{1,2,\ldots,2m\}$. The boundary edges of this blocker contain
the edges $[j,j+1]$ and $[i_t,i_{t}+1]$, but not the edges
$[i_t-1,i_t]$ and $[2m-1,2m]$. Hence, they are not consecutive on
the boundary $\langle 1,2,\ldots,2m,1 \rangle$ of this subgraph,
contradicting the inductive assumption.

\textbf{Case 4: $E=\{[0,1],[i_t,i_t+1]\}$.} If $i_t=1$, we are
done. If not, then $B \setminus \{[0,1]\}$ is a blocker in the
subgraph spanned by $\{1,2,\ldots,2m\}$, whose boundary edges are
$[i_t,i_t+1]$, and possibly $[2m,1]$. If both belong to $B
\setminus \{[0,1]\}$, this contradicts the inductive assumption,
since $i_t \neq 1$, and thus the boundary edges are not
consecutive. Otherwise, $B \setminus \{[0,1]\}$ contains only one
boundary edge of the induced subgraph, and this is impossible for
a blocker (see Figure~\ref{fig:T7}).
\end{proof}

\subsection{Characterization of the Blockers}

In this section we present a complete characterization of the
blockers and prove one direction of Theorem~\ref{Thm:Main}. We
start with a few observations.

\begin{figure}[tb]
\begin{center}
\scalebox{0.6}{
\includegraphics{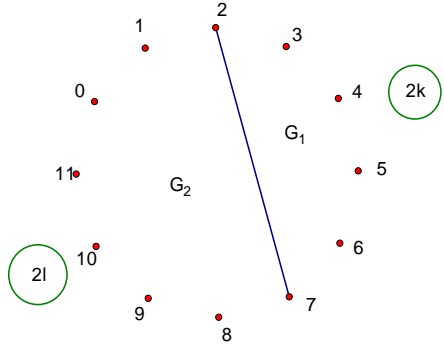}
} \caption{An edge in a blocker separates the remaining vertices
into two sets. In this figure, $m=6$, $e=[a_i,a_j]=[2,7]$, $k=2$
and $l=3$.} \label{fig:T8}
\end{center}
\end{figure}

Let $B$ be a blocker in $CK(2m)$, and let $e=[i,j] \in B$ ($0 \leq
i < j \leq 2m-1$) be an edge that separates the remaining vertices
into sets of $2k$ and $2l$ vertices, respectively. Clearly, $k,l
\geq 0$ and $k+l=m-1$ (see Figure~\ref{fig:T8}).

Denote by $G_1^{-}$ the subgraph of $CK(2m)$ of order $2k$ spanned
by the vertices $\{i+1,i+2,\ldots,j-1\}$, and by $G_1^{+}$ the
subgraph of $CK(2m)$ of order $2k+2$ spanned by the vertices
$\{i,i+1,i+2,\ldots,j-1,j\}$. Similarly, denote by $G_2^{-}$ the
subgraph of $CK(2m)$ of order $2l$ spanned by the vertices
$\{j+1,j+2,\ldots,2m-1,0,1,\ldots,i-1\}$, and by $G_2^{+}$ the
subgraph of $CK(2m)$ of order $2l+2$ spanned by the vertices
$\{j,j+1,j+2,\ldots,2m-1,0,1,\ldots,i-1,i\}$.

We observe that in $G_1^{-}$, there exists an SPM all whose edges
are parallel to $e$: $\{[i+\nu,j-\nu] | \nu=1,2,\ldots,k\}$.
Hence, $G_2^{+}$ does not include an SPM disjoint from $B$, since
otherwise the union of these two SPMs would form an SPM in
$CK(2m)$ that misses $B$. In other words, this implies that $B
\cap E(G_2^{+})$ is a blocking set in $G_2^{+}$. In particular, it
follows that $|B \cap E(G_2^{+})| \geq l+1$, and thus $|(B
\setminus \{e\}) \cap E(G_2^{+})| \geq l$.

Repeating the same argument with $G_1^{-}$ replaced by $G_2^{-}$
and $G_2^{+}$ replaced by $G_1^{+}$, we find that $B \cap
E(G_1^{+})$ is a blocking set in $G_1^{+}$, hence $|B \cap
E(G_1^{+})| \geq k+1$ and $|(B \setminus \{e\}) \cap E(G_1^{+})|
\geq k$.

Since $|B|=m=k+l+1$ and the edge sets $(B \setminus \{e\}) \cap
E(G_1^{+}), (B \setminus \{e\}) \cap E(G_2^{+})$, and $\{e\}$ are
pairwise disjoint, it follows that
\[
|(B \setminus \{e\}) \cap E(G_1^{+})|= k, \qquad \qquad |(B
\setminus \{e\}) \cap E(G_2^{+})|=l,
\]
and that $B$ does not contain an edge that crosses $e$ in an
interior point. Since this argument holds for any $e \in B$, we
conclude that $B$ is a {\it simple} (i.e., crossing-free) set of
edges of $CK(2m)$.

Now we are ready to prove the main theorem. First, we recall its
statement:
\begin{theorem}\label{Thm:Main'}
Any blocker in $CK(2m)$ is a simple caterpillar graph whose spine
lies on the boundary of $CK(2m)$ and is of length $t \geq 2$. If
the spine ``starts'' with the vertex $0$ and the edge $[0,1]$,
then the edges of the blocker are:
\[
\{[i-1,i]:1 \leq i \leq t\} \cup
\{[t+j-1-\epsilon_{t+j},t+j+\epsilon_{t+j}]:1 \leq j \leq m-t\},
\]
where the $\epsilon_i$ are integers that satisfy $1 \leq
\epsilon_{t+1} < \epsilon_{t+2} < \ldots < \epsilon_m \leq m-2$.
\end{theorem}
\begin{proof}
Let $B$ be a blocker in $CK(2m)$. Denote the number of its
boundary edges by $t$, and assume that these edges ``start'' with
the vertex $0$ and the edge $[0,1]$. By Lemma~\ref{Lemma:Main},
the boundary edges of $B$ are $\{[0,1],[1,2],\ldots,[t-1,t]\}$. We
make the following three observations:
\begin{enumerate}

\item $B$ cannot contain an edge of the form $[l_1,l_2]$ for $1
\leq l_1+1<l_2 \leq t$, since such an edge is either of even
order, or is parallel to one of the boundary edges of $B$.

\item $B$ cannot contain a non-boundary edge $e$ such that all
boundary edges of $B$ lie on one side of $e$. Indeed, if this
happens, then we can define the subgraphs $G_1^{+}$ and $G_2^{+}$
as in the beginning of this subsection, and find that $B \cap
E(G_1^{+})$ is a blocker in $G_1^{+}$ that contains only one
boundary edge of $G_1^{+}$, which is impossible, since $|G_1^{+}|
\geq 4$.

\medskip

A combination of these two observations implies that any
non-boundary edge of $B$ connects one of the vertices
$1,2,\ldots,t-1$ with one of the vertices $t+1,t+2,\ldots,2m-1$.
Furthermore, the only edge in $B$ that contains the vertex $0$ is
$[0,1]$, and the only edge in $B$ that contains the vertex $t$ is
$[t-1,t]$.

\medskip

\item Consider two non-boundary edges $[l_1,l'_1],[l_2,l'_2] \in
B$, and assume that $0 < l_1 < l_2 < t$. (Thus, by the previous
observation, $t < l'_1,l'_2 < 2m$). First note that $l'_2 \leq
l'_1$, since the edges of $B$ do not cross. Secondly we show that
$l'_1 - l'_2 > l_2-l_1$. Indeed, if the differences are equal then
the edges $[l_1,l'_1]$ and $[l_2,l'_2]$ are parallel, which is
impossible for a blocker. If $l'_1 - l'_2 < l_2-l_1$, then one can
construct an SPM that misses $B$. The edges of the SPM are the
following (see Figure~\ref{fig:T9}):

\begin{enumerate}

\item All the edges parallel to $[l_2,l'_2]$ on its left side:
$[l_2-\epsilon,l'_2+ \epsilon]$, for $1 \leq \epsilon \leq
m-\frac{1}{2}(l'_2-l_2+1)$. (Addition and subtraction here are
modulo $2m$.)

\item All the edges parallel to $[l_1,l'_1]$ that lie on the right
side of the edge $[l_2,l'_2]$ (in the weak sense): $[l_2 +
\epsilon, l'_1-l_2+l_1-\epsilon]$, for $0 \leq \epsilon \leq
(l_1+l'_1-1)/2 - l_2$.

\item Alternating boundary edges (this is possible since the
number of remaining vertices is $l'_2 - (l'_1-l_2+l_1) =
(l'_2+l_2) - (l'_1+l_1)$, which is a positive even number).
\end{enumerate}
\end{enumerate}
The assertion of the theorem follows immediately from these three
observations.

\begin{figure}[tb]
\begin{center}
\scalebox{0.6}{
\includegraphics{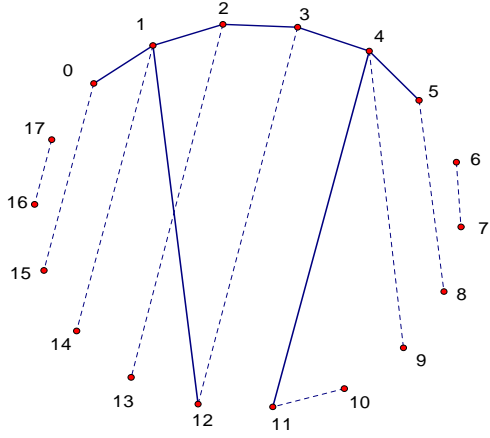}
} \caption{An illustration of Observation~3 in the proof of
Theorem~\ref{Thm:Main'}. In this figure, $m=9$, $t=5$, $l_1=1$,
$l_2=4$, $l'_1=12$, and $l_2'=11$. The SPM that misses $B$ is
depicted by dotted lines.} \label{fig:T9}
\end{center}
\end{figure}

The exact formulation of the theorem is explained as follows. The
theorem lists the $m$ edges of $B$. The first $t$ are the boundary
edges, and the remaining $m-t$ are the non-boundary edges,
arranged by decreasing order of the vertex in which they meet the
``boundary path'' of $B$, $\{[0,1],\ldots,[t-1,t]\}$. Denoting by
$[l_2,l'_2]$ the $(t+j)$-th edge in the list, and by $[l_1,l'_1]$
the $(t+j+1)$-st edge, we get:
\[
l_1 = t+j-\epsilon_{t+j+1}, \qquad l'_1 = t+j+1+\epsilon_{t+j+1},
\]
\[
l_2 = t+j-1-\epsilon_{t+j}, \qquad l'_2 = t+j+\epsilon_{t+j},
\]
and therefore,
\[
l'_1-l'_2 = 1+\epsilon_{t+j+1}-\epsilon_{t+j} >
\epsilon_{t+j+1}-\epsilon_{t+j}-1 = l_2-l_1.
\]
\end{proof}

\subsection{Proof of the Inverse Direction}
\label{sec:sub:inverse}

In this subsection we prove the inverse direction of
Theorem~\ref{Thm:Main}, namely, that any caterpillar subgraph of
$CK(2m)$ that satisfies the conditions mentioned above is a
blocker.

%This definition will help us to determine whether two diagonals of
%$P_{2m}$ that belong to some sub-polygon $P_{2l}$ are ``parallel''
%in that polygon.

%\medskip

%\noindent Now we are ready to prove the inverse direction of
%Theorem~\ref{Thm:Main}.
\begin{theorem}\label{Thm:Main''}
Let $B$ be the following set of $m$ edges of $CK(2m)$:
\[
\{[i-1,i]:1 \leq i \leq t\} \cup
\{[t+j-1-\epsilon_{t+j},t+j+\epsilon_{t+j}]:1 \leq j \leq m-t\},
\]
where the $\epsilon_i$ are integers that satisfy $1 \leq
\epsilon_{t+1} < \epsilon_{t+2} < \ldots < \epsilon_m \leq m-2$
(and hence, $t \geq 2$). Then $B$ is a blocker in $CK(2m)$.
\end{theorem}
\begin{proof}
The proof uses double induction: A primary induction on $m$, and
a secondary (backward) induction on the number of boundary edges
in $B$.

\medskip

\noindent For $m=2$ the claim is clear. For $m=3$ there are only
two possible sets $B$ that satisfy the conditions (up to
isomorphism). The first is a path of three consecutive boundary
edges, which is indeed a blocker in $CK(6)$ by
Observation~\ref{obs:basic}. The second consists of all diagonals
of odd order emanating from a single vertex, and therefore
intersects every SPM in one edge.

\medskip

\noindent For $m \geq 4$, we assume that the assertion holds for
$m-1$ and prove it for $m$. Let $B$ be a set of $m$ edges of
$CK(2m)$ satisfying the assumptions, and let $t$ be the number of
boundary edges in $B$. If $t=m$, then $B$ is a path of $m$
consecutive boundary edges of $CK(2m)$, which is indeed a blocker
by Observation~\ref{obs:basic}. If $t<m$, we assume that the
assertion holds for all sets $B$ satisfying the assumptions and
having more than $t$ boundary edges, and prove the assertion for
$B$. Assume on the contrary that $B$ is not a blocker, and thus
there exists an SPM $M$ that does not meet $B$. We distinguish two
cases:\footnote{Recall that by the assumptions of
Theorem~\ref{Thm:Main''}, the boundary edges of $B$ are
$[0,1],[1,2],\ldots,[t-1,t]$.}
\begin{enumerate}
\item \textbf{Case A: $[2m-1,0] \in M$.} In this case we omit the
vertices $2m-1$ and $0$ from $CK(2m)$, and show that $B'=B
\setminus \{[0,1]\}$ satisfies the assumptions of the theorem for
the induced subgraph $CK(2m-2)$ spanned by the vertices
$1,2,\ldots,2m-2$. On the other hand, $B'$ is not a blocker, since
$M'=M \setminus \{[2m-1,0]\}$ is an SPM in that graph that does
not meet $B'$. This contradicts the inductive assumption on $m$.

\item \textbf{Case B: $[2m-1,0] \not \in M$.} In this case we add
to $B$ the edge $[2m-1,0]$, and omit from $B$ the edge parallel to
$[2m-1,0]$. We obtain a new set $B''$ that also satisfies the
assumptions of the theorem (for the same graph $CK(2m)$), and has
$t+1$ boundary edges. On the other hand, $B''$ is not a blocker,
since it does not meet $M$, contradicting the inductive assumption
on $t$.
\end{enumerate}
Now we discuss the two cases in more detail:
\begin{enumerate}
\item \textbf{Case A.} If $t=2$ then $B$ is the set of all edges
of odd order emanating from the vertex $1$. In this case, $B'$ is
the set of all edges of odd order in $CK(2m-2)$ emanating from the
vertex $1$, and thus clearly satisfies the assumption of
Theorem~\ref{Thm:Main''}.

If $t>2$ then $B'$ is a caterpillar whose spine lies on the
boundary of $CK(2m-2)$. The spine contains the edges
$[1,2],[2,3],\ldots,[t-1,t]$. Each of the other edges of $B$,
i.e., the edges of the form
$\{[t+j-1-\epsilon_{t+j},t+j+\epsilon_{t+j}]:1 \leq j \leq m-t\}$,
is ``parallel'' to the respective boundary edge $[t+j-1,t+j]$ also
in $CK(2m-2)$. Thus, $B'$ contains a representative of each of the
$m-1$ directions in $CK(2m-2)$.

If $\epsilon_m=m-2$ then the last edge of the form
$[t+j-1-\epsilon_{t+j},t+j+\epsilon_{t+j}]$ is $[1,2m-2]$, which
is a boundary edge in $CK(2m-2)$. This edge extends the path
$\langle 1,2,\ldots, t \rangle$ from the left, and thus the length
of the spine of $B'$ is $t$. If $\epsilon_m<m-2$, then all the
non-boundary edges of $B$ are also non-boundary edges in
$CK(2m-2)$.

In addition, any non-boundary edge of $B'$ (with respect to
$CK(2m-2)$) connects one of the internal vertices of the spine
with one of the internal vertices of the rest of the boundary of
$CK(2m-2)$, and if $e_1=[p_1,q_1]$ and $e_2=[p_2,q_2]$ are two
non-boundary edges of $B'$, where $p_1,p_2$ are on the spine and
$q_1,q_2$ are not on the spine, then the distance between $q_1$
and $q_2$ is greater than the distance between $p_1$ and $p_2$
also with respect to $CK(2m-2)$ (i.e., this property is also
inherited from the properties of $B$ in $CK(2m)$). This shows that
$B'$ satisfies the assumptions of Theorem~\ref{Thm:Main''} with
respect to $CK(2m-2)$.

\item \textbf{Case B.} It is clear from the construction that
$B''$ contains a representative of each of the $m$ directions in
$CK(2m)$. In addition, $B''$ is a caterpillar whose spine $\langle
2m-1,0,1,\ldots,t \rangle$ is a path on the boundary of $CK(2m)$.
We have to check that the non-boundary edges of $B''$ connect an
internal vertex of the spine with an internal vertex of the rest
of the boundary of $CK(2m)$. In order to verify this, it is
sufficient to show that $B$ does not contain an edge of the form
$[2m-1,i]$ for $0 < i \leq t$. This is indeed true since such an
edge is either of even order or parallel to one of the boundary
edges of $B$, hence cannot be contained in $B$. The other
condition checked in Case A (the comparison of the distances
between $q_1,q_2$ and $p_1,p_2$ for two non-boundary edges
$[p_1,q_1]$ and $[p_2,q_2]$) holds for $B$, hence also for $B''$.
This shows that $B''$ satisfies the assumptions of
Theorem~\ref{Thm:Main''} with respect to $CK(2m)$.
\end{enumerate}

This completes the proof of Theorem~\ref{Thm:Main''}, and with it
the proof of Theorem~\ref{Thm:Main}.
\end{proof}

\subsection{The Number of Blockers in $CK(2m)$}

\begin{proposition}
The number of blockers in $CK(2m)$ is $m \cdot 2^{m-1}$.
\end{proposition}
\begin{proof}
We partition the set of blockers in $CK(2m)$ into subsets
according to the number $t$ of boundary edges, and count the
number of blockers with exactly $t$ boundary edges. Consider the
blockers whose spine (i.e., path of boundary edges) is $\langle
0,1,\ldots,t \rangle$. Each such blocker is uniquely determined by
the number of edges emanating from each of the vertices
$1,2,\ldots,t-1$. Hence, the number of such blockers is equal to
the number of nonnegative integer solutions of the equation
$x_1+x_2+\ldots+x_{t-1}=m-t$, that is known to be
${{m-t+t-2}\choose{t-2}}={{m-2}\choose{t-2}}$. Thus, the number of
blockers whose spine is of the form $\langle 0,1,\ldots \rangle$
is $\sum_{t=2}^{m} {{m-2}\choose{t-2}} = 2^{m-2}$. Therefore, by
symmetry, the total number of blockers in $CK(2m)$ is $2m \cdot
2^{m-2} = m \cdot 2^{m-1}$, as asserted.
\end{proof}

\end{document}